\documentclass{article}
\usepackage[english]{babel}
\usepackage[centertags]{amsmath}
\usepackage{amsfonts,amssymb,amsthm}
\usepackage{float,longtable}
\usepackage[latin1]{inputenc}
\usepackage{newlfont}

\newcommand{\Om}{\Omega}
\newcommand{\la}{\langle}
\newcommand{\ra}{\rangle}

\newenvironment{pf}{\noindent{\sc Proof}.\enspace}{\rule{2mm}{2mm}\medskip}

\newtheorem{theorem}{Theorem}
\newtheorem{proposition}{Proposition}
\newtheorem{lemma}{Lemma}
\newtheorem{corollary}{Corollary}
\newtheorem{remark}{Remark}
\newtheorem{remarks}{Remark}
\newtheorem{definition}{Definition}
\newcommand{\be}{\begin{equation}}
\newcommand{\ee}{\end{equation}}
\newcommand{\teta}{\vartheta}
\renewcommand{\theta}{\vartheta}
\newcommand{\om}{\omega}

\newcommand{\ov}{\bar}

\newcommand{\R}{\mathbb{R}}   

\newcommand{\N}{\mathbb{N}}   
\newcommand{\T}{\mathbb{T}}   

\renewcommand{\a }{\alpha }
\renewcommand{\b }{\beta }
\newcommand{\s }{\sigma }

\renewcommand{\d }{\delta }

\newcommand{\e }{\varepsilon }
\newcommand{\g }{\gamma}

\renewcommand{\l }{\lambda }

\newcommand{\vphi}{\varphi }

\newcommand{\ph}{\varphi}

\newcommand{\ub}{\bar{u}}   
\newcommand{\vb}{\bar{v}}

\newcommand{\sn}{\mathrm{sn}}  
\newcommand{\dn}{\mathrm{dn}}  
\newcommand{\am}{\mathrm{am}}  

\newcommand{\Vb}{\bar{V}}
\newcommand{\Ob}{\bar{\Omega}}
\newcommand{\mb}{\bar{m}}

\newcommand{\dps}{\displaystyle}

\long\def\symbolfootnote[#1]#2{\begingroup%
\def\thefootnote{\fnsymbol{footnote}}\footnote[#1]{#2}\endgroup}   

\begin{document}

\title{\textbf{
Periodic 
solutions of 
wave equations for asymptotically 
full measure sets of frequencies}}

\author{Pietro Baldi\footnote{Sissa, via Beirut 2-4, 34014, Trieste, Italy.\, 
\textit{E-mail:} \texttt{baldi@sissa.it}\,.},
\,  Massimiliano Berti\footnote{Dipartimento di Matematica e Applicazioni ``R. Caccioppoli'', Universit\`a di Napoli ``Federico II'' , via Cintia, 80126, Napoli, Italy.\,
\textit{E-mail:} \texttt{m.berti@unina.it}\,.} 
}
\date{}

\maketitle

%
%

\section{Introduction}

The aim of this Note is to prove existence and multiplicity 
of small amplitude periodic solutions of
the completely resonant wave 
equation\symbolfootnote[0]{\emph{Keywords:} Nonlinear Wave Equation, 
Infinite dimensional Hamiltonian Systems, 
Periodic solutions, Lyapunov-Schmidt reduction, Small divisors problem.

\emph{2000AMS Subject Classification:} 35L05, 35B10, 37K50.

Supported by MURST within the PRIN 2004 
``Variational methods and nonlinear differential equations''.}  
\be\label{eq:main}
\begin{cases} \square u + f(x,u) =  0 \cr
u ( t, 0 )= u( t,  \pi )  = 0
\end{cases}
\ee
where $ \square := \partial_{tt} - \partial_{xx} $
is the D'Alambertian operator and 
\be\label{nonlin}  
f(x,u) = a_2 u^2 + a_3(x)u^3 + O( u^4 ) 
\qquad \mathrm{or} \qquad  \ 
f(x,u) = a_4 u^4 + O( u^5 )
\ee
for a Cantor-like set of frequencies $\om $
of asymptotically full measure at 
$\om = 1$. 

Equation (\ref{eq:main}) is called completely resonant because
any solution $v = \sum_{j \geq 1} a_j \cos ( j t + \teta_j ) \sin (jx ) $
of the linearized equation at $ u = 0 $
\be\label{eq:lin}
\begin{cases} 
u_{tt} - u_{xx} =  0 \cr
u ( t, 0 )= u( t,  \pi )  = 0
\end{cases}
\ee
is $2 \pi$-periodic in time.
\\[1mm]
\indent
Existence and multiplicity of periodic solutions of completely resonant
wave equations had been proved 
for a zero measure, uncontable Cantor set of frequencies in
\cite{BP1} for $ f(u) = u^3 + O(u^5 )$ and in 
\cite{BB}-\cite{BB1} for
any nonlinearity $ f(u) = a_p u^p +O(u^{p+1})$, $ p \geq 2 $. 

Existence of periodic solutions for a Cantor-like set 
of frequencies of asymptotically full measure  
has been recently proved in \cite{BB2}  where, 
due to the well known ``small divisor difficulty",
the ``$0$th order bifurcation equation'' is required to possess 
non-degenerate periodic solutions. Such property was
verified in 
\cite{BB2} for nonlinearities 
like $ f = a_2 u^2 + O(u^4) $, 
$  f  = a_3(x) u^3 + O(u^4 ) $. 
See also \cite{GMP} for 
$f= u^3 + O(u^5 )$.
\\[1mm]
\indent
In this Note we shall prove that, for quadratic, cubic and quartic 
nonlinearities $ f(x,u) $ like in (\ref{nonlin}),
the corresponding $0$th order bifurcation equation possesses non-degenerate 
periodic solutions -- Propositions \ref{63} and 
\ref{prop caso u2u3} --,
implying, by the results of  \cite{BB2},    Theorem \ref{thm:main} and Corollary \ref{FR} below. 

We remark that our proof is purely analytic (it does not use
numerical calculations) being based on the analysis of the variational equation
and exploiting properties of the Jacobi elliptic functions.

\subsection{Main results}

Normalizing the period to $ 2 \pi $, we look for solutions 
of  
\begin{equation*}  
\begin{cases}
\om^2 u_{tt} - u_{xx} + f (x, u ) = 0 \cr
u(t,0)= u(t, \pi) = 0
\end{cases}
\end{equation*}
in the Hilbert algebra 
(for $ s > 1 \slash 2 $, $ \s > 0  $) 
\begin{eqnarray*}
X_{\s,s} := \Big\{  u (t,x) = \sum_{l \geq 0} 
\cos (lt) \ u_l (x)   & \Big| &  u_l  \in H^1_0 ((0,\pi), \R) \ \ 
 \forall l \in \N \ \  \mathrm{and} \\ 
& &  
\|u\|_{\s,s}^2 := \sum_{l \geq  0} \exp{(2\s l)} 
(l^{2s} + 1 ) \|u_l\|^2_{H^1} < +\infty \Big\} \, .
\end{eqnarray*}

It is natural to look for solutions which are even in time
because equation (\ref{eq:main}) is reversible. 

We look as well for solutions of (\ref{eq:main})  
in the subalgebras
\begin{eqnarray*}
X_{\s,s,n} := \Big\{  u \in X_{\s,s} \ | \ u \ \text{is} \ 
\frac{2\pi}{n}\text{-periodic} \Big\} \subset X_{\s,s} \, , \quad n \in \N 
\end{eqnarray*}
(they are particular $2\pi$-periodic solutions). 

The space of the solutions of the linear equation (\ref{eq:lin}) 
that belong to 
$H^1_0( \mathbb{T} \times (0, \pi), \R ) $ and  are even in time
is
\begin{eqnarray*}
V &:= &  \Big\{  v (t,x) = \sum_{l \geq 1}  \cos (lt) 
u_l \sin (lx) \  \Big|  \ u_l \in \R \, , \    
\sum_{l \geq 1} l^2 | u_l |^2 < +\infty \Big\} \\
& = & \Big\{ v(t,x) = \eta ( t + x ) - \eta ( t - x ) \ 
 \Big| \  \eta  \in H^1(\T, \R) 
 \ \mathrm{with} \  \eta \ \mathrm{odd} \Big\} \, .
\end{eqnarray*}

\begin{theorem}\label{thm:main} 
Let
\be\label{fu2u3}
f(x,u) = 
a_2 u^2 + a_3(x) u^3 + \sum_{k\geq 4} a_k(x) u^k 
\ee
where 
$(a_2, \la a_3 \ra ) \neq (0,0)$, 
$\la a_3 \ra := \pi^{-1}  
\int_0^\pi a_3 (x) dx $, or
\be\label{fu4}
f(x,u) = 
a_4 u^4 + \sum_{k\geq 5} 
a_k(x) u^k  
\ee
where $ a_4 \neq 0 $, $ a_5 ( \pi - x ) = - a_5(x) $, 
$ a_6( \pi - x) = a_6(x) $, $ a_7( \pi - x) = - a_7(x) $. 
Assume moreover  $ a_k (x) \in H^1((0, \pi), \R) $ with
 $ \sum_{k} \| a_k \|_{H^1}  \rho^k < + \infty $ for some $ \rho > 0 $. 

Then there exists $ n_0 \in \N$ such that $ \forall n \geq n_0 $ there is 
$ \d_0 > 0 $, 
$ \ov{\s} > 0 $ and a $C^\infty$-curve 
$[0, \d_0) \ni \delta \to u_{\delta} \in X_{{\ov \sigma} \slash 2,s,n}$
with the following properties:
\begin{itemize}
\item $(i)$ 
$ \left\| u_{\delta} - \delta {\ov v}_n \right\|_{{\ov \s}\slash 2, s,n} = 
O( \delta^2 )$ 
for some 
$ {\ov v}_n \in V \cap X_{{\ov \s}, s,n} \setminus \{0\}$ 
with minimal period $2 \pi \slash n $;
\item  $(ii)$ 
there exists  a Cantor set ${\cal C}_n 
\subset [0,\delta_0)$ of asymptotically full measure, i.e. 
satisfying 
\be\label{meas}
\lim_{\varepsilon \to 0^{+}} \frac{\mathrm{meas} ({\cal C}_n 
\cap (0, \varepsilon))}{\varepsilon} = 1 \, ,
\ee
such that, $\forall $ $ \delta \in
{\cal C}_n $, 
$ u_{\delta} ( \om (\d) t,x) $ is a $2\pi/ (\om(\d) n) $-periodic,
classical solution of (\ref{eq:main}) with 
$$
\om (\d)= 
\begin{cases}  \sqrt{1-2 s^* \d^2} \qquad \text{if $f$ is like in 
(\ref{fu2u3})} 
\cr 
\sqrt{1-2 \d^6} \qquad \ \ \, \text{if $f$ is like in (\ref{fu4})}
\end{cases}
$$ 
and\footnote{Note how the interaction between the second and the third 
order terms $a_2 u^2 $, 
$ a_3 (x) u^3 $ 
changes the bifurcation diagram, i.e.
existence of periodic solutions for frequencies $\om $ less or/and greater  
of $ \om = 1$.} 
$$
s^* = \begin{cases}
 -1 \qquad \mathrm{if} \qquad \la a_3 \ra \geq \pi^2 a_2^2 \slash 12  
\cr \pm 1 
\qquad \mathrm{if} \qquad 0 < \la a_3 \ra < \pi^2 a_2^2 \slash 12 \cr 
1 \ \ \qquad \mathrm{if} \qquad \la a_3 \ra \leq 0 \, .
\end{cases}
$$
\end{itemize}
\noindent 
\end{theorem}

By (\ref{meas}) also each Cantor-like 
set of frequencies 
$ {\cal W}_n := \{ \om (\d) \ | \ \d \in {\cal C}_n \}$ has asymptotically
full measure at $ \om = 1 $.

\begin{corollary}\label{FR}    {(\bf Multiplicity)}
There exists a Cantor-like set 
${\cal W}$ of asymptotically full measure at $ \om = 1 $, such that
$ \forall \om \in {\cal C} $, equation (\ref{eq:main}) 
possesses geometrically distinct periodic solutions
$$
u_{n_0}  \, , \ldots \, , 
u_n \,  , \ldots  \,    u_{N_\om} \, , \qquad N_\om \in \N 
$$
with the same period $ 2\pi/ \om $.
Their number increases arbitrarily as $ \om $ tends to $ 1 $:
 $$
 \lim_{\om \to 1}  N _ \om  =  + \infty  \, . 
 $$
\end{corollary}

\begin{pf}
The proof is like in \cite{BB2} and we report it for completeness.
If $ \d $ belongs to 
the  asymptotically full measure set (by (\ref{meas}))
$$
D_n := {\mathcal C}_{n_0} \cap 
\ldots \cap  {\mathcal C}_n \, ,  \qquad n \geq n_0 
$$
there exist $(n - n_0+1)$ geometrically distinct 
periodic solutions of (\ref{eq:main}) with the same 
period $ 2\pi/ \om (\d) $ 
(each $u_n $ has minimal period $ 2\pi/ (n\om (\d)) $). 

There exists  a decreasing sequence  of positive $ \e_n \to 0 $ such that
$$
 \mathrm{meas} (D_n^c \cap (0 ,  \e_n ) ) \leq  \e_n 2^{-n}  \,  .
 $$
Let define the set $ {\mathcal C }  \equiv D_n $ 
on each $[ \e_ {n+1}, \e_n)$. 
$  {\mathcal C } $ has asymptotically full measure at $ \d = 0 $ 
and for each $ \d \in  {\mathcal C }  $
there exist $  N ( \d ) := \max   \{  n  \in \N :     \d  <  \e_n \} $ geometrically distinct 
periodic solutions of  (\ref{eq:main})  with the 
same period  $ 2\pi/ \om (\d)$.  $ N(\d) \to +\infty $ as $\d \to 0 $. 
\end{pf}

\begin{remark}
Corollary \ref{FR} is an analogue for equation (\ref{eq:main}) of 
the well known multiplicity results of Weinstein-Moser \cite{We}-\cite{Mo}
and Fadell-Rabinowitz  \cite{FR} 
which hold in finite dimension. 
The 
solutions  form a sequence of  functions 
with increasing norms and decreasing minimal periods. 
Multiplicity of solutions 
was also obtained in \cite{BB1}  
(with the "optimal" number $ N_\omega \approx C \slash  \sqrt{|\om - 1 |}$) 
 but only for a zero measure set of 
frequencies. 
\end{remark}

The main point for proving Theorem \ref{thm:main} relies in showing
the existence of non-degenerate solutions 
of the 
$0$th order bifurcation equation for $ f $ like in (\ref{nonlin}).
In these cases the $0$th order bifurcation equation 
involves higher order terms of the nonlinearity, 
and, for $n$ large, can be reduced to an integro-differential equation
(which physically describes 
an averaged effect  of the nonlinearity with
Dirichlet boundary conditions).
\\[2mm]
{\bf Case $f(x,u) = a_4 u^4 + O(u^5 )  $.}  Performing the rescaling
$$
u \to \d u \, , \qquad \d > 0 
$$
we look for $2\pi \slash n $-periodic solutions  in $   X_{\s,s,n}  $ of 
\be\label{eq:freqre}
\begin{cases}
\om^2 u_{tt} - u_{xx} + \d^3 g (\d, x, u ) = 0 \cr
u(t,0)= u(t, \pi) = 0
\end{cases}
\ee
where 
$$
g (\d, x, u) :=  \frac{f(x,\d u) }{\d^4} = 
a_4  u^4 + \d a_{5} (x) u^{5} + \d^2 a_{6} (x) u^{6}
+ \ldots \, .
$$
To find solutions of (\ref{eq:freqre})
we implement the Lyapunov-Schmidt reduction 
according to the orthogonal decomposition 
$$
X_{\s,s,n} = (V_n \cap X_{\s,s,n}) \oplus (W \cap X_{\s,s,n})
$$ 
where
\begin{eqnarray*}
V_n & := 
\Big\{ v(t,x) = \eta ( n t + n x ) - \eta ( nt -  nx ) \ 
| \  \eta  \in H^1(\T, \R) \ \  \mathrm{with} \ \ 
\eta \ \mathrm{odd} \Big\}  
\end{eqnarray*}
and
\begin{equation*}  
W :=
\Big\{ w = \sum_{l\geq 0 } \cos (lt)\  w_l(x) \in X_{0,s} \ | \ 
\ \ 
\int_0^\pi w_l (x) \sin (lx) \, dx = 0, \ \forall l \geq  0 \ \Big\} \, .
\end{equation*}
Looking for solutions 
$ u = v + w $ with $ v \in  V_n \cap X_{\s,s,n} $, $ w \in W  \cap X_{\s,s,n} $,
and imposing the frequency-amplitude relation 
$$
\frac{(\om^2 - 1 )}{2} = - \d^6 
$$
we are led to solve the bifurcation equation 
and the range equation 
\be\label{eqs1} 
\begin{cases}
\Delta v = \d^{-3} 
\Pi_{V_n}  g (\d, x,v + w)  
\cr 
\dps L_\om w = \d^3 \Pi_{W_n}  g(\d, x, v + w) 
\ \qquad  
\end{cases}
\ee
where 
$$ 
 \Delta v := v_{xx} + v_{tt}, \qquad 
\qquad   
L_\om := - \om^2 \partial_{tt} + \partial_{xx} 
$$ 
and 
$ \Pi_{V_n}  : X_{\s,s,n} \to V_n \cap X_{\s,s,n}  $,  $ \Pi_{W_n}  : X_{\s,s,n} \to W  \cap X_{\s,s,n}  $
denote the projectors. 

With the further rescaling 
$$
w \to \d^3 w
$$
and since 
$ v^4 \in W_n $ (Lemma 3.4 of \cite{BB}), 
$ a_5 (x) v^5 $, $ a_6 (x) v^6 $, $ a_7 (x) v^7 \in W_n $
because
$ a_5 ( \pi - x ) = - a_5(x) $, 
$ a_6( \pi - x) = a_6(x) $, $ a_7( \pi - x) = - a_7(x) $
(Lemma 7.1 of \cite{BB2}), system (\ref{eqs1}) 
is equivalent to
\be \label{nPQu2}
\left\{ \begin{array}{ll}
\Delta v = \Pi_{V_n}   \Big( 4 a_4 v^3 w + \d 
r(\d,x, v, w ) \Big)  \\
L_\om w = a_4 v^4 + 
\d \Pi_{W_n}  
 \widetilde r (\d,x, v, w) 
\end{array} \right.
\ee
where $r(\d, x,v,w) = 
a_8(x) v^8+ 5 a_5(x) v^4 w + O(\d )$ and 
$ \widetilde r(\d, x,v,w) = a_5 (x) v^5 + O(\d)  $. 

For $ \d = 0 $ system (\ref{nPQu2}) reduces to 
$ w = - a_4 \square^{-1} v^4 $ and to  
the $0$th order bifurcation equation 
\be\label{0bif}
\Delta v + 4 a_4^2 \Pi_{V_n}   \Big( v^3 \square^{-1} v^4 \Big)  = 0  
\ee
which is  
the Euler-Lagrange equation of
the functional  $\Phi_0 : V_n  \to \R$ 
\be \label{Phi0u4}
\Phi_0 (v) = \frac{\|v\|_{H_1}^2}{2} - \frac{a_4^2}{2} 
\int_\Om v^4 \square^{-1} v^4 
\ee
where $\Omega := \T \times (0, \pi)$. 

\begin{proposition}\label{63} Let $ a_4 \neq 0 $.  
$ \exists $ $ n_0 \in \N $ such that $ \forall n \geq n_0 $ 
the $0$th order bifurcation equation 
(\ref{0bif}) has a  solution $ \ov{v}_n  \in V_n  $
which is non-degenerate in $ V_n $ (i.e. $ \mathrm {Ker}  D^2 \Phi_0 = \{ 0 \}$), 
with minimal period $ 2 \pi \slash n $.
\end{proposition}

\noindent
{\bf Case $f(x,u) = a_2 u^2 + a_3(x)u^3 + O( u^4 ) . $} 
Performing the rescaling
$u \to \d u$ 
we look for $2\pi \slash n $-periodic solutions  of 
\begin{equation*}  
\begin{cases}
\om^2 u_{tt} - u_{xx} + \d g (\d, x, u ) = 0 \cr
u(t,0)= u(t, \pi) = 0
\end{cases}
\end{equation*}
where 
$$
g (\d, x, u) :=  \frac{f(x,\d u) }{\d^2} = 
a_2  u^2 + \d a_3 (x) u^3 + \d^2 u_4(x) u^4 + \dots  \, . 
$$
With the frequency-amplitude relation 
$$
\frac{\om^2-1}{2} = - s^* \d^2
$$
where $ s^* = \pm 1 $, we have to solve
\be \label{PQu2}
\left\{  \begin{array}{l}
- \Delta v= - s^* \d^{-1} \Pi_{V_n}  g(\d , x, v + w )   \\
L_\om w= \d \Pi_{W_n}  g(\d, x, v+w) \,. 
\end{array} \right. 
\ee
With the further rescaling $w \to \d w $\,
and since $ v^2 \in W_n  $, system (\ref{PQu2}) 
is equivalent to
\be \label{nPQu2 seconda volta} 
\Bigg\{ \begin{array}{l}
-\Delta v = s^* 
\Pi_{V_n}  \Big( -2a_2 v w - a_2 \d  w^2 - a_3(x) (v+ \d w)^3 
- \d r(\d,x, v+\d w) \Big) \\
L_\om w = a_2 v^2 + \d \Pi_{W_n}  \Big( 2a_2 v w + \d a_2 w^2+ 
a_3 (x) (v+ \d w )^3 + \a_8(x) v^8d r(\d,x, v+\d
w) \Big) 
\end{array} 
\ee
where $r(\d, x,u):=\d^{-4} [f(x,\d u)-a_2\d^2u^2 - \d^3 
a_3(x) u^3]=a_4(x) u^4 + \dots $

For $ \d = 0 $ system (\ref{nPQu2 seconda volta}) reduces to
$ w = -a_2 \square^{-1} v^2 $ and the 0th order bifurcation equation
\be\label{0th order bif eq caso u2u3}
- s^* \Delta v
= 2a_2^2 \Pi_{V_n}  (  v \square^{-1} v^2) - \Pi_{V_n}  (a_3 (x) v^3) 
\ee
which is the Euler-Lagrange equation of  $\Phi_0 : V_n  \to \R$ 
\begin{equation}\label{phiu2u3} 
\Phi_0 (v) := 
s^* \frac{  \| v \|_{H^1}^2}{2} - \frac{a_2^2}{2} \int_\Om v^2 \square^{-1} v^2
+ \frac{1}{4} \int_\Om a_3(x) v^4  \, . 
\end{equation}

\begin{proposition}\label{prop caso u2u3} Let $ (a_2, \la a_3 \ra ) \neq 0 $. 
$ \exists $ $ n_0 \in \N $ such that $ \forall n \geq n_0 $ 
the $0$th order bifurcation equation 
(\ref{0th order bif eq caso u2u3}) has a  solution $ \ov{v}_n  \in V_n  $
which is non-degenerate in $ V_n $, with minimal period $ 2 \pi \slash n $.
\end{proposition}

\section{Case $ f ( x , u ) = a_4 u^4 + O( u^5 )$}

We have to prove the existence of {\it non-degenerate} critical points 
of
the functional
$$ 
\Phi_n : V  \to \R  \, , \qquad  
\Phi_n (v) := \Phi_0 ({\cal H}_n v ) 
$$ 
where $\Phi_0 $ is defined  in (\ref{Phi0u4}). Let
$ {\cal H}_n : V \to V $ be the linear isomorphism 
defined, for $ v(t,x) = \eta (t+x) - \eta (t-x) \in V $,
by 
$$a_8(x) v^8
({\cal H}_n v)(t,x) := \eta (n (t+x)) - \eta (n (t-x)) 
$$
so that $V_n \equiv {\cal H}_n V $.

\begin{lemma} {\bf See \cite{BB1}.}
$\Phi_n $ has the following development: 
for $ v(t,x) = \eta (t+x) - \eta (t-x) \in V $ 
\be\label{phin}
\Phi_n ( \b n^{1\slash 3} v ) = 
4 \pi \b^2 n^{8\slash 3}   
\Big[ \Psi (\eta) + \a \frac{ {\cal R}(\eta)}{n^2} \Big] 
\ee
where $ \b := (3 \slash (\pi^{2} a_4^{2}))^{1\slash 6 }$, $ \a := a_4^2 \slash (8\pi) $, 
\be\label{psi}
\Psi (\eta ) := 
\frac{1}{2} \int_{\T} \eta'^2 (t) \, dt - 
\frac{2\pi}{8} \Big(  \la \eta^4 \ra + 3 \la \eta^2 \ra^2 \Big)^2 \, , 
\ee
$ \la \  \ra $ denotes the average on $ \T $, and 
\be\label{Re}
{\cal R}(\eta ) := - \int_\Om v^4 \square^{-1} v^4 \, dt dx + 
\frac{\pi^4}{6} 4 \Big(  \la \eta^4 \ra + 3 \la \eta^2 \ra^2 \Big)^2 \, .
\ee
\end{lemma}

\begin{pf}
Firstly the quadratic term writes 
\be\label{kin}
\frac{1}{2} \| {\cal H}_n v \|_{H^1}^2 = \frac{n^2}{2} 
\| v \|_{H^1}^2 =
 n^2 2 \pi \int_{\T} \eta'^2 (t) \, dt \, . 
\ee
By Lemma 4.8 in \cite{BB1} 
the non-quadratic term can be developed as
\be\label{lemn}
\int_\Om ({\cal H}_n v)^4 \square^{-1} ({\cal H}_n v)^4 = 
\frac{\pi^4}{6}  \la m \ra^2
- \frac{{\cal R}(\eta )}{n^2} 
\ee
where $ m : \T^2 \to \R  $ is 
$ m (s_1, s_2 ) := $ $ ( \eta (s_1) - \eta (s_2)  )^4 \, $,
$ \la m \ra := (2\pi)^{-2} \int_{\T^2} m(s_1, s_2) \, ds_1 d s_2 $ 
denotes its average, and 
\be\label{Reta}  
{\cal R}(\eta ) :=  \Big( - \int_\Om v^4 \square^{-1} v^4 + 
\frac{\pi^4}{6} \la m \ra^2 \Big) 
\ee
is homogeneous of degree $ 8 $.
Since $ \eta $ is  odd
we find 
\be\label{ave}
\la m \ra = 2 \Big(  \la \eta^4 \ra + 3 \la \eta^2 \ra^2 \Big) 
\ee
where $\la \  \ra $ denotes the average on $ \T$.

Collecting (\ref{kin}), (\ref{lemn}), (\ref{Reta}) and (\ref{ave}) we find out 
\begin{equation*}    
\Phi_n ( \eta ) = 2 \pi n^2 \int_{\T} \eta'^2 (t) \, dt \, - 
\frac{\pi^4}{3} a_4^2
\Big(  \la \eta^4 \ra + 3 \la \eta^2 \ra^2 \Big)^2  + 
\frac{a_4^2}{2n^2}{\cal R}(\eta ) \, .
\end{equation*}
Via the rescaling $ \eta \to \beta n^{1\slash 3} \eta $
we get the expressions (\ref{psi}) and (\ref{Re}).
\end{pf}

By (\ref{phin}), in order 
to find for $ n $ 
large enough a non-degenerate critical
point of $ \Phi_n $, 
it is sufficient to find a non-degenerate critical point of
$\Psi (\eta ) $ defined 
on 
$$
E := \Big\{ \eta \in H^1(\T), \, \eta \ \mathrm{odd} \Big\} \, , 
$$
namely non-degenerate solutions in $ E $ of 
\be\label{eq:u4}
{\ddot \eta} + A(\eta) 
\Big( 3 \la \eta^2 \ra \eta + \eta^3 \Big) = 0\, \ \ \qquad
A(\eta) := \la \eta^4 \ra + 3 \la \eta^2 \ra^2  \, .
\ee

\begin{proposition}
There exists an odd,  analytic, $2\pi$-periodic solution $ g(t)  $ 
of (\ref{eq:u4}) which is non-degenerate in $ E $. 
$ g(t) = V \sn (\Om t, m)$\, 
where $\sn$ is the Jacobi elliptic sine and  
\, $V >0 $,\, $\Om>0$, $m \in (-1,0)$ are suitable constants
(therefore $g(t)$ has minimal period $2 \pi $).
\end{proposition}

We will construct the solution $ g $ of (\ref{eq:u4}) 
by means of the Jacobi elliptic sine  in Lemma \ref{jac}. 
The existence of a solution $ g $ 
follows also directly applying to $ \Psi : E \to \R  $ 
the Mountain-Pass Theorem \cite{AR}. Furthermore such solution is an 
analytic function arguing as in Lemma 2.1 of \cite{BB2}. 

\subsection{Non-degeneracy of $g$}

We now want to prove that $ g $ is non-degenerate. The linearized 
equation of (\ref{eq:u4}) at $g $ is
\begin{eqnarray*}    
& & {\ddot h} + 3A( g)\Big[ \la g^2 \ra h + g^2 h \Big] +
6 A(g)g \la gh\ra + A'(g)[h] \Big( 3 \la g^2 \ra g + g^3 \Big) 
  =  \\
& & {\ddot h} + 3A( g)\,\big[ \la g^2 \ra  + g^2  \big]\,h +
6g \la gh\ra \Big(\la g^4 \ra + 3 \la g^2 \ra^2 \Big) +
4g \Big( \la g^3h\ra + 3\la g^2\ra \la gh \ra \Big)
\Big( 3 \la g^2 \ra + g^2 \Big) 
  =  0 
\end{eqnarray*}
that we write as
\be\label{linea0}
{\ddot h} + 3A( g)\, \big( \la g^2 \ra  + g^2  \big)\, h = - \la gh \ra I_1 - 
\la g^3 h\ra I_2 
\ee
where
\begin{eqnarray}   
\begin{array}{l}   \label{I1}
\begin{cases}
 I_1 := 6\Big( 9 \la g^2 \ra^2 + \la g^4 \ra \Big) g +12 \la g^2\ra g^3 \vspace{4pt} \\
 I_2 := 12 g \la g^2 \ra + 4 g^3 \, .
\end{cases}
\end{array}
\end{eqnarray}
For $ f \in E $, let $ H := L(f) $ be the unique solution belonging to $ E $  
of the non-homogeneous linear system
\be\label{inve}
{\ddot H} + 3A( g)\,\big( \la g^2 \ra  + g^2  \big)\,H = f \, ;
\ee
an integral representation of the Green operator 
$L$ is given in Lemma \ref{lemma della formula integrale di L}. Thus
(\ref{linea0}) becomes
\be\label{eql}
 h = -\la gh \ra L( I_1) -\la g^3 h \ra L(I_2) \, .
\ee
Multiplying (\ref{eql}) by $g$ and taking averages we get
\begin{eqnarray}  \label{incorniciata 1}
	\la gh \ra \,\big[ 1+ \la gL(I_1) \ra \big] = - \la g^3 h \ra\,\la gL(I_2)\ra, 
\end{eqnarray}
while multiplying (\ref{eql}) by $g^3$ and taking averages
\begin{eqnarray}   \label{incorniciata 2}
	\la g^3h \ra \big[ 1+ \la g^3 L(I_2) \ra \big] = - \la gh \ra\,\la g^3 L(I_1)\ra.
\end{eqnarray}
Since $g$ solves (\ref{eq:u4}) we have the following identities.

\begin{lemma}
There holds
\be\label{id1}
2 A(g) \la g^3 L(g) \ra = \la g^2 \ra \,
\ee
\be\label{id2}
2 A(g) \la g^3 L( g^3) \ra = \la g^4 \ra \, .
\ee
\end{lemma}

\begin{pf}
(\ref{id1}) 
is obtained by the identity for $ L(g) $ 
$$
\frac{d^2}{dt^2}(L(g)) + 3A(g)\,\big( \la g^2 \ra  + g^2  \big)\,L(g) = g 
$$
multiplying by $ g $, taking averages, integrating by parts, 
$$
\la {\ddot g} L(g) \ra 
+ 3A(g)\,\big[ \la g^2 \ra \la L(g)g \ra + \la g^3 L(g) \ra \big] = \la g^2 \ra 
$$
and using that $ g $ solves (\ref{eq:u4}).

Analogously, (\ref{id2}) is obtained by the identity for $ L(g^3) $
$$
\frac{d^2}{dt^2}(L(g^3)) + 3A(g)\,\big( \la g^2 \ra  + g^2  \big)\,L(g^3)= g^3 
$$
multiplying by $ g $, taking averages, integrating by parts, 
 and using that $ g $ solves (\ref{eq:u4}).
\end{pf}

Since $L$ is a symmetric operator 
we can compute the following averages using (\ref{I1}), (\ref{id1}), (\ref{id2}):
\begin{eqnarray}
\begin{cases}
  \la g L(I_1) \ra = 6 \,\Big( \la g^4 \ra +9 \la g^2 \ra^2 \Big)\, 
\la gL(g) \ra + 6\,A(g)^{-1}\, \la g^2 \ra^2 \vspace{2pt}\\   \label{coefficienti}
  \la g L(I_2) \ra = 12 \la g^2 \ra \,\la gL(g) \ra + 2\,A(g)^{-1}\, \la g^2 \ra 
\vspace{3pt} \\
  \la g^3 L(I_1) \ra = 9 \la g^2 \ra \vspace{4pt} \\
  \la g^3 L(I_2) \ra = 2 \, .
\end{cases}
\end{eqnarray}
Thanks to the identities  
(\ref{coefficienti}), equations (\ref{incorniciata 1}), (\ref{incorniciata 2}) 
simplify to
\begin{eqnarray}   
\begin{cases}
\la gh \ra \, 
\big[ A(g)+6\la g^2 \ra^2 \big] B(g) = -2\,\la g^2 \ra\, B(g)\,\la g^3 h \ra  \\
\la g^3 h \ra = -3 \la g^2 \ra\,\la gh \ra  \label{nuova incorniciata 2}
\end{cases}
\end{eqnarray}
where 
\begin{eqnarray}
	B(g):= 1+ 6 A(g) \la gL(g) \ra \, .
\end{eqnarray}
Solving (\ref{nuova incorniciata 2}) we get 
\[ B(g) \la gh \ra =0 \, .
\]
We will prove in Lemma \ref{glg}\, that $B(g) \neq 0$, so $\la gh \ra =0$. Hence by (\ref{nuova incorniciata 2}) also $\la g^3 h \ra =0$\, and therefore, by (\ref{eql}),\, $h=0$.  This concludes the proof of the non-degeneracy of the solution $g$ of (\ref{eq:u4}).

\vspace{9pt}
It remains to prove that $B(g) \neq0$.  The key is to express the function 
$ L(g) $ by means of the variation of constants formula.

We first look for a fundamental set of solutions of the homogeneous equation
$$
{\ddot h} + 3A(g) \,\big( \la g^2 \ra + g^2 \big)\, h = 0 \, . \eqno{\mathrm{(HOM)}}
$$ 

\begin{lemma}   \label{uv}
There exist two linearly independent solutions  of (HOM),
${\ov u} := {\dot g}(t) \slash {\dot g}(0)$ and ${\ov v}$, such that
\[
\begin{cases} 
{\ov u} \ \ \mathrm{is \  even, \ 2\pi \, periodic}  \\ 
{\ov u} (0) = 1 \, , \ {\dot {\ov u}}(0) = 0 
\end{cases} 
\qquad 
\begin{cases} 
{\ov v} \ \ \mathrm{is \ odd, \: not \: periodic}\\ 
{\ov v} (0) = 0 \, , \ {\dot {\ov v}}(0) = 1 
\end{cases}
\]
and
\be \label{rhopos}
\vb (t+2\pi )-\vb (t) = \rho \ub (t) \, \qquad \mathrm{for \ some \ } \rho > 0 \, .
\ee
\end{lemma}

\begin{pf}
Since (\ref{eq:u4}) is autonomous, $ {\dot g}(t) $ is a solution of 
the linearized equation (HOM). $ {\dot g}(t) $ is even and $ 2 \pi $-periodic. 

We can construct another solution of (HOM) in the following way. 
The super-quadratic Hamiltonian system 
(with constant coefficients)
\begin{eqnarray}   \label{eq coeff cost}
\ddot{y}+3A(g)\la g^2 \ra \, y + A(g)\, y^3=0  
\end{eqnarray}
possesses a one-parameter family of 
odd, $T(E)$-periodic solutions $ y (E, t ) $,\,  close to $ g $, 
parametrized by the energy $ E$. Let ${\ov E} $ denote the energy level of  
$ g $,\, i.e. $g = y ( \ov E, t )$\, and $T(\ov E ) = 2\pi  $.
 
Therefore 
$ l(t) := ( \partial_E y(E,t) )_{ |E = {\ov E}}$ is an odd solution of (HOM). 

 Deriving the identity $ y (E, t+T(E)) = y(E, t )  $\, 
with respect to $E$ we obtain at $ E=\bar{E} $
\[
l(t+2 \pi) - l(t) = - (\partial_E T(E))_{|E={\ov E}}\, {\dot g}(t)
\]
and, normalizing $ \ov{v}(t) := l(t) / \dot{l}(0)$, we get (\ref{rhopos})
with
\be\label{rhoint}
\rho := -(\partial_E T(E))_{|E={\ov E}}\,
\Big( \frac{\dot{g}(0)}{\dot{l}(0)} \Big) \, .
\ee
Since $ y(E,0) =0 $ $ \forall E $, the energy identity gives
$ E = \frac{1}{2} ({\dot y}(E,0))^2 $. Deriving
w.r.t $ E $ at $E = {\ov E}$, yields 
$ 1= {\dot g}(0) {\dot l}(0) $ which, inserted in (\ref{rhoint}), gives
\be\label{TE}
\rho = -(\partial_E T(E))_{|E={\ov E}}\, (\dot{g}(0))^2 \, . 
\ee
$\rho > 0 $ because $ (\partial_E T(E))_{|E={\ov E}} < 0  $ by the
superquadraticity of the potential of (\ref{eq coeff cost}). 
It can be checked also by a computation, see Remark after Lemma \ref{jac}.
\end{pf}

Now we write an integral formula for the Green operator $ L $.

\begin{lemma}  \label{lemma della formula integrale di L}
For every $f \in E $ there exists a 
unique solution $ H = L ( f )$ of (\ref{inve}) which can be written as 
\be  \label{Formula di L}
L (f) = 
\Big( \int_0^t f(s) {\ov u}(s)\, ds + \frac{1}{\rho } \int_0^{2\pi} f {\ov v}  
\Big) {\ov v} (t) -
\Big( \int_0^t f(s) {\ov v}(s)\, ds \Big) {\ov u} (t) \ \in E \, .
\ee
\end{lemma}

\begin{pf}
The non-homogeneous equation (\ref{inve})
possesses the particular solution
$$
{\ov H}(t) = 
\Big( \int_0^t f(s) {\ov u}(s)\, ds \Big) {\ov v} (t) -
\Big( \int_0^t f(s) {\ov v}(s)\, ds \Big) {\ov u} (t) 
$$
as can be verified noting that
the Wronskian $ {\ov u}(t){\dot {\ov v}} (t) 
- {\dot {\ov u}}(t){\ov v} (t) \equiv 1 $, $ \forall t $.
Notice that $\ov H $ is odd. 

Any solution $H(t)$\, of (\ref{inve}) can be written as
$$
H(t)={\ov H}(t) + a {\ov u} + b {\ov v} \, , \qquad a, b \in \R \, .
$$
Since $\ov H $ is odd, ${\ov u}$ is even and ${\ov v}$ is odd,  
requiring $ H $ to be odd, implies 
$ a = 0 $.
Imposing now the $ 2 \pi $-periodicity yields
\begin{eqnarray}
0 \hspace{-7pt} & = \hspace{-8pt} & \Big( \int_0^{t + 2\pi} f \ov u \Big) {\ov v }(t+2 \pi ) -
\Big( \int_0^{t + 2\pi} f \ov v \Big) {\ov u }(t+2 \pi ) - 
\Big( \int_0^{t} f \ov u \Big) {\ov v }(t) +
\Big( \int_0^t f \ov v \Big) {\ov u }(t) +
b \Big( {\ov v}(t+2\pi) - {\ov v}(t) \Big) \nonumber \\
& = \hspace{-8pt} & \Big( b + \int_0^t f \ov u \Big) 
\Big( {\ov v }(t+2 \pi ) - {\ov v }(t) \Big) - {\ov u }(t)
\Big( \int_t^{t + 2\pi} f \ov v \Big)                         \label{ult}
\end{eqnarray}
using that $ \ov u $ and $ f \ov u $ are $2\pi$-periodic and 
$ \la f \ov u \ra = 0 $.
By (\ref{ult}) and (\ref{rhopos}) we get
\be\label{f1}
\rho \Big( b + \int_0^t f \ov u \Big)  - \int_t^{t + 2\pi} f \ov v   = 0\, .
\ee
The left hand side in (\ref{f1}) is constant in time
because, deriving w.r.t. $t$,
$$
\rho f(t) \ov u (t) -  f(t) \Big( \ov v(t+{2\pi}) - \ov v (t) \Big) = 0  
$$
again by (\ref{rhopos}).\, 
Hence evaluating (\ref{f1})\, for $t= 0 $\, yields
\,$ b = $ $\rho^{-1}  \int_0^{2\pi} f {\ov v} $.\, So there exists a 
unique solution $ H = L ( f )$ of (\ref{inve}) belonging to $E$ and
(\ref{Formula di L})\, follows.
\end{pf}

Finally

\begin{lemma}There holds  \label{glg}
$$
\la g L(g) \ra = \frac{\rho }{4 \pi A(g)}  +
\frac{1}{2 \pi \rho } \Big( \int_0^{2\pi} g {\ov v} \Big)^2 > 0
$$
because $ A (g )$, $ \rho > 0 $.
\end{lemma}

\begin{pf}
Using (\ref{Formula di L}) we can compute 

\begin{eqnarray}
\la g L(g) \ra & = &  
\frac{1}{2\pi} \int_0^{2\pi} 
\Big( \int_0^t g {\ov u} \Big) {\ov v}(t) g(t) \, dt \ + \ 
\frac{1}{2\pi \rho} \Big( \int_0^{2\pi} g {\ov v} \Big)^2 \ - \frac{1}{2\pi} \int_0^{2\pi}
\Big( \int_0^t g{\ov v}\Big) {\ov u} (t) g(t) \, dt \nonumber \\
& = &
2 \frac{1}{2\pi} \int_0^{2\pi}
\Big( \int_0^t g {\ov u} \Big) {\ov v}(t) g(t) \, dt \  +
\frac{1}{2\pi \rho} \Big( \int_0^{2\pi} g {\ov v} \Big)^2 \label{finq}
\end{eqnarray}
because, by $ \int_0^{2\pi} g \ov u = 0 $, we have
$$
0 = \int_0^{2\pi}
\frac{d}{dt} \Big[ \Big( \int_0^t g {\ov v} \Big)\Big( \int_0^t g {\ov u} \Big)\, \Big] dt \   
= \int_0^{2\pi}
\Big[ \Big( \int_0^t g {\ov v} \Big) {\ov u}(t) g(t)  +
 \Big( \int_0^t g {\ov u} \Big) {\ov v}(t) g(t) \Big] \,  dt  \, .
$$
Now, since $ \ov u (t)= {\dot g}(t) \slash {\dot g}(0) $ and $ g(0) = 0 $,
\be\label{pe}
\frac{1}{2\pi} \int_0^{2\pi}
\Big( \int_0^t g {\ov u} \Big) {\ov v}(t) g(t)  = \frac{1}{2\pi {\dot g}(0)}
\int_0^{2\pi} 
\Big( \int_0^t \frac{d}{d\tau} \frac{g^2(\tau)}{2}\, d\tau \Big) {\ov v}(t) g(t)  =
\frac{1}{4\pi{\dot g}(0)}  \int_0^{2\pi} g^3  {\ov v}  \, .
\ee
We claim that
\be\label{last}
\int_0^{2\pi} g^3  {\ov v}  = \frac{\rho {\dot g}(0)}{2 A(g)}\, .
\ee
By (\ref{finq}), (\ref{pe}), (\ref{last}) we have the thesys.

Let us prove (\ref{last}). Since $ g $ solves (\ref{eq:u4}) multiplying
by $ \ov v $ and integrating
\be\label{s1}
\int_0^{2\pi} \ov v(t) {\ddot g}(t) 
+ 3 A(g) \la g^2 \ra g(t) \ov v(t) 
+ A(g) g^3(t) \ov v (t) \, dt  = 0
\ee
Next, 
since $\ov v $ solves (HOM), multiplying
by $ g $ and integrating
\be\label{s2}
\int_0^{2\pi} g(t) {\ddot {\ov v}}(t) + 3 A(g) \la g^2 \ra 
{\ov v}(t)g(t) + 3 A(g) g^3(t) {\ov v}(t) \, dt = 0 \, .
\ee
Subtracting (\ref{s1}) and (\ref{s2}), gives 
\be\label{r1}
\int_0^{2\pi} \ov v(t) {\ddot g}(t) -  g(t) {\ddot {\ov v}}(t) 
= 2 A(g) \int_0^{2\pi} g^3  {\ov v}  \, .
\ee
Integrating by parts the left hand side, since $g(0)= g(2\pi ) = 0 $,\, $\ub (0)=1$\, 
and (\ref{rhopos}),   gives
\be\label{r2}
\int_0^{2\pi} \ov v(t) {\ddot g}(t) -  g(t) {\ddot {\ov v}}(t) 
= {\dot g}(0) 
 [ v(2\pi ) - v(0) ] = \rho {\dot g} (0)\, .
\ee
(\ref{r1}) and (\ref{r2}) give (\ref{last}).
\end{pf}

\subsection{Explicit computations}

We now give the explicit construction of $g$ by means of 
the Jacobi elliptic sine defined as follows.
Let $\am ( \cdot, m) : \R \to \R $ be the inverse function of the 
Jacobi elliptic integral of the first kind
$$
\vphi \mapsto F (\vphi, m ) := 
\int_{0}^\vphi \frac{d\theta}{\sqrt{1-m \sin^2 \! \theta}}\, .
$$
The Jacobi elliptic sine is defined by 
$$
\sn (t, m) := \sin ( \am (t, m))\, .
$$
$ \sn(t,m) $ is $ 4 K(m) $-periodic, where $K(m)$ is the complete 
elliptic integral of the first kind
\[
K(m) :=F \Big( \frac{\pi}{2}, m \Big) 
= \int_0^{\pi /2} \frac{d\theta}{\sqrt{1-m \sin^2\! \theta }}
\]
and admits an analytic extension with a pole in
$iK( 1 - m ) $ for $m \in (0,1) $ and in 
$iK\big( 1/(1-m) \big) /  \sqrt{1-m}$ for $m < 0 $.
\, Moreover, since 
\[
\partial_t \am (t,m)=\sqrt{ 1-m\, \sn^2(t,m)}\,,
\]
the elliptic sine satisfies
\begin{equation}   \label{sn}
(\dot{\sn})^2=(1-\sn^2)(1-m\,\sn^2)\,.
\end{equation}
  
\begin{lemma}  \label{jac}
There exist 
\,$ V >0$,\, $\Om>0 $, 
$ m \in (-1,0)$\, 
such that $ g(t):= V \sn (\Om t, m)$\, 
is an odd, analytic, $2\pi$-periodic solution of (\ref{eq:u4})
with pole in 
$iK\big( 1/(1-m) \big) /\big( \Omega \,\sqrt{1-m} \big)$.
\end{lemma}

\begin{pf} 
Deriving (\ref{sn}) we have 
\,$\ddot{\sn} +(1+m)\,\sn  -2m\, \sn^3  = 0 $.\,
Therefore $ g_{(V,\Om,m)}(t) := V \sn (\Omega t , m )$\,
is an odd, $(4K(m)/\Om)$-periodic solution of
\begin{equation}  \label{eq di Vsn(Omega t,m)} 
{\ddot g} + \Omega^2 (1+m) g -  2 m \frac{\Omega^2}{V^2} g^3 = 0 \, .
\end{equation}
The function $g_{(V,\Om,m)} $ will be a solution of (\ref{eq:u4}) if
$(V,\Om,m)$ verify 
\begin{eqnarray}   \label{ipotesi H}
\begin{cases}
  	\Omega^2 (1+m) =  3 A(g_{(V,\Om,m)})\,\la g_{(V,\Om,m)}^2 \ra \cr
   -2m \Om^2 = V^2 A(g_{(V,\Om,m)}) \cr
   2 K(m)= \Om \pi \, .  
   \end{cases}
\end{eqnarray}
Dividing the first equation of (\ref{ipotesi H}) 
by the second one
\be\label{eq:m}
-\frac{1+m}{6m}= \la \sn^2 (\cdot, m) \ra \, .
\ee
The right hand side can be expressed as  
\be  \label{media di sn2}
\la \sn^2 (\cdot, m) \ra = \frac{K(m)-E(m)}{mK(m)}
\ee
where $ E(m)$ is the complete elliptic integral of the second kind
$$
E(m) :=  \int_0^{\pi /2}  \sqrt{1-m \sin^2\! \theta }\,d\theta
= \int_0^{K(m)} 1 - m \, \sn^2(\xi,m) \, d \xi  
$$
(in the last passage we
make the change of variable $\theta = \am (\xi, m) $).

Now, we show that system (\ref{ipotesi H}) has a 
unique solution.
By (\ref{eq:m}) and (\ref{media di sn2}) 
\begin{equation}  \label{eq per m}
(7+m)K(m) -6E(m)=0 \, .
\end{equation}
By the definitions of $ E(m)$ and $ K(m) $ we have
\be\label{psim}
\psi(m) := (7+m)K(m) -6E(m) 
=  \int_0^{\pi\slash 2} \frac{1+m(1+6\sin^2\! \theta)}
{\big( 1-m\sin^2 \!\theta \big)^{1/2}} \, d \theta \, .
\ee
For $m=0$ it holds \,$ \psi(0) 
= \pi /2 >0$ and, for $ m = - 1 $, 
$ \psi(-1) = 
- \int_0^{\pi /2} 6 \sin^2\! \theta \, (1 + \sin^2\! \theta)^{-1\slash 2} \, d\theta <0 .$
Since $\psi $ is continuous there exists a solution $ \mb \in (-1,0) $ 
of (\ref{eq per m}). 
Next the third equation in (\ref{ipotesi H}) fix $ \Ob $ and 
finally we  find $ \Vb $. Hence 
$ g(t) =\Vb\,\sn(\Ob t,\mb)$ solves (\ref{eq:u4}).

Analyticity and poles follow from \cite{Handbook}, 16.2, 16.10.2, pp.570,573.

At last, $\mb $ is unique because $\psi'(m) > 0 $  
for $ m \in (-1,0) $ as can be verified by (\ref{psim}).
One can also compute that $\mb \in (-0.30, -0.28)$.
\end{pf}


\noindent
\textbf{Remark.}\ We can compute explicitly the sign of $dT/dE$ and $\rho$ of (\ref{TE}) in the following way.

The functions $ g_{(V,\Om,m)} $ are 
solutions of 
the Hamiltonian system (\ref{eq coeff cost}) 
imposing 
\begin{equation}   \label{ipotesi H_0}
\begin{cases}	
 \Om^2 (1+m)= \alpha   \\
  -2 m \Om^2 = V^2 \beta 
  \end{cases}
\end{equation}
where $\alpha := 3 A(g)\,\la g^2 \ra $,\,  $\beta := A(g)$ and 
$ g $ is the solution constructed in Lemma \ref{jac}. 

We solve (\ref{ipotesi H_0}) w.r.t $ m $ finding 
the one-parameter family $(y_m)$ of odd 
periodic solutions $ y_m(t) := V(m)\, \sn (\Om(m) t ,\, m ) $,\, close to $g$,  with 
energy and period
$$
E(m) = \frac{1}{2} V^2(m) \Om^2(m)\, =\, - \frac{1}{\b}\, m\, \Om^4(m)\, ,  \qquad \quad
T(m) = \frac{4K(m)}{\Om (m)} \, . 
$$ 
It holds
$$
\frac{dT(m)}{dm} = \frac{4K'(m)\Om(m) - 4 K(m) \Om' (m)}{\Om^2(m)} > 0 
$$
because $K'(m) > 0 $\, and from (\ref{ipotesi H_0}) 
\,$\Om'(m) = - \Om (m) \big( 2 (1+m) \big)^{-1} < 0 $.\, Then    
$$
\frac{dE(m)}{dm} = - \frac{1}{\b}\, \Om^4(m) - \frac{1}{\b} \,m\, 4 \Om^3(m) \Om'(m) < 0\, ,
$$
so
\[
\frac{dT}{d E} = \frac{d T(m)}{dm} \Big( \frac{d E(m)}{dm} \Big)^{-1} < 0 \,
\]
as stated by general arguments in the proof of Lemma \ref{uv}.

We can also write an explicit formula for $\rho$,
\be\label{rho finale}
\rho = \frac{m}{m -1} \bigg[ 2\pi 
+ (1+m) \int_0^{2\pi} \frac{\sn^2 (\Om t,m)}{\dn^2 (\Om t,m)}\, d t \bigg]\,. 
\ee
From (\ref{rho finale}) it follows that $\rho >0$\, because \,$-1 < m <0 $.

\section{
Case $f(x,u) = a_2 u^2 + a_3(x) u^3 +O( u^4 )$}

We have to prove the existence of {\it non-degenerate} critical points of the functional
$ \Phi_n (v) := \Phi_0 ({\cal H}_n v )  $ where $ \Phi_0 $ is defined in 
(\ref{phiu2u3}).
 
\begin{lemma} \label{lemma: decomp di Phi0 con Psi caso u2u3}
See \cite{BB1}. $\Phi_n$ has the following development: for $v(t,x) = \eta (t+x) - \eta (t-x) \in V$,
\begin{equation} \label{eq decomp Phi caso u2u3}
\Phi_n ( \b n v) = 4 \pi \beta^2
 n^4 \Big[ \Psi(\eta) +
\frac{\b^2}{4\pi} 
\Big( 
\frac{R_2(\eta)}{ n^2} + R_3 (\eta) \Big) \Big]
\end{equation}
where 
\[
\Psi(\eta):=\frac{s^* }{2}  \int_{{\T}} {\dot \eta}^2 + \frac{\beta^2}{4 \pi} \Big[ \a \Big(\int_{{\T}} \eta^2 \Big)^2 + \g  \int_{{\T}} \eta^4 \Big]
\]
\begin{equation}  \label{R2R3}
R_2 (\eta )  := 
-\frac{a_2^2}{2} \Big[ \int_\Om v^2 \square^{-1} v^2 - \frac{\pi^2}{6}
\Big( \int_{{\T}} \eta^2  \Big)^2 \Big],
\quad \quad
R_3 (\eta) :=
\frac{1}{4} \int_\Om \big( a_3(x) - \la a_3 \ra \big) ({\cal H}_n v)^4 \,, 
\end{equation}
$\alpha:=\big( 9\la a_3 \ra - \pi^2 a_2^2 \big)/12$,\: 
$\gamma:=\pi \la a_3 \ra /2$,  and
\[
\beta=
\left\{
\begin{array}{ll}
\!\! (2|\alpha|)^{-1/2} & \mathrm{if}\;\, \a \neq 0, \\
\!\! (\pi/\g)^{1/2}     & \mathrm{if}\;\, \a = 0.
\end{array}     \right.
\]
\end{lemma}

\begin{pf}
By Lemma 4.8 in \cite{BB1} with 
$m(s_1,s_2)=\big( \eta(s_1)-\eta(s_2) \big)^2$, 
for $v(t,x) = \eta (t+x) - \eta (t-x)$ the operator  $ \Phi_n $ admits the development
\begin{eqnarray*}
\Phi_n (v)  &= &  2\pi s^* n^2 \int_{{\T}} {\dot \eta}^2 (t) dt 
- \frac{\pi^2 a_2^2}{12} \Big( \int_{{\T}} \eta^2 (t) \, dt \Big)^2 -
\frac{a_2^2}{2n^2}  \Big( \int_\Om v^2 \square^{-1} v^2 - \frac{\pi^2}{6}
\Big( \int_{{\T}} \eta^2 (t) \, dt \Big)^2 \, \Big) \\
&  & + \,\frac{1}{4} \la a_3 \ra \int_\Om v^4 + 
\frac{1}{4} \int_\Om \big( a_3(x) - \la a_3 \ra \big) ({\cal H}_n v)^4\, .
\end{eqnarray*}
Since
$$
\int_\Om v^4 = 2 \pi \int_{{\T}} \eta^4  + 3 
\Big( \int_{{\T}} \eta^2 \Big)^2  ,
$$
we write
\begin{eqnarray*}
\Phi_n (v)  &= &  
2\pi s^* n^2 \int_{{\T}} {\dot \eta}^2 (t) dt 
- \frac{\pi^2 a_2^2}{12} \Big( \int_{{\T}} \eta^2 \Big)^2 +
\frac{1}{4} \la a_3 \ra \Big[ 2 \pi \int_{{\T}} \eta^4  + 3 
\Big( \int_{{\T}} \eta^2 \Big)^2 \Big] +
\frac{R_2(\eta)}{n^2} + R_3 (\eta)\, , 
\end{eqnarray*}
where $R_2$, $R_3$ defined in (\ref{R2R3}) are both homogenous of degree 4. 
So
\begin{eqnarray*}
\Phi_n (v) = 2\pi s^* n^2 \int_{{\T}} {\dot \eta}^2 + \a
\Big( \int_{{\T}} \eta^2 \Big)^2 + 
\g \int_{{\T}} \eta^4  +
\frac{R_2(\eta)}{n^2} + R_3 (\eta)
\end{eqnarray*}
where $\a$, $\g$ are defined above.
With the rescaling $\eta \to \eta \b n $ we get 
decomposition (\ref{eq decomp Phi caso u2u3}). \end{pf}

In order to find for $n$ large a non-degenerate critical point of $\Phi_n$, 
by (\ref{eq decomp Phi caso u2u3}) it is sufficient 
to find critical points of $\Psi$\, on 
$E = \big\{ \eta \in H^1(\T), \, \eta \ \mathrm{odd} \big\}$ 
(like in Lemma 6.2 of \cite{BB2} also the term $ R_3(\eta) $
tends to $ 0 $ with its derivatives). 

If $\la a_3 \ra \in (-\infty,0) \cup (\pi^2 a_2^2 /9, +\infty)$, then $\a \neq 0$ and 
we must choose  $s^*=-\mathrm{sign}(\a)$, so that the functional becomes
\begin{equation*}
\Psi (\eta)= \mathrm{sign}(\a) \Big( -\frac12 \int_{\T} \dot{\eta}^2 + \frac{1}{8\pi} \Big[ \Big( \int_{\T} \eta^2 \Big)^2 + \frac{\g}{\a} \int_{\T} \eta^4 \Big] \Big)\,.
\end{equation*}
Since in this case $\g/\a >0$, 
the functional $\Psi$ clearly has a mountain pass critical point, 
solution of
\begin{equation} \label{ODE per valori esterni}
\ddot{\eta}+\la \eta^2 \ra \eta + 
\lambda \eta^3 =0\,, \qquad \qquad \lambda=\frac{\g}{2\pi \a} >0\,.
\end{equation}
The proof of the non-degeneracy of the solution of 
(\ref{ODE per valori esterni}) is very simple using the analytical arguments
of the previous section (since $\l > 0 $ it is sufficient a positivity argument).


If $\la a_3 \ra =0$, then the equation becomes 
\,$\ddot{\eta}+\la \eta^2 \ra \eta =0$,\, so we find 
again what proved in \cite{BB2} for $ a_3(x) \equiv 0 $.

If $\la a_3 \ra =\pi^2 a_2^2 /9$,\, then $\a=0$. 
We must choose $s^*=-1$, so that we obtain
\[
\Psi(\eta)=-\frac12 \int_{\T} \dot{\eta}^2 + \frac14 \int_{\T} \eta^4,
\qquad \quad \ddot{\eta}+\eta^3=0.
\]
This equation has periodic solutions which are non-degenerate because of  
non-isocronicity, see Proposition 2 in \cite{BP}.

Finally, if $ \la a_3 \ra  \in (0,\, \pi^2 a_2^2 /9)$,\, 
then $\a<0$ and there are both solutions for 
$ s^* = \pm 1$. The functional
\begin{eqnarray*}
\Psi (\eta) &= &
\frac{s^*}{2} \int_{\T} \dot{\eta}^2 + 
\frac{1}{8\pi} \Big[ -\Big( \int_{\T} \eta^2 \Big)^2 + 
\frac{\g}{|\a|} \int_{\T} \eta^4 \Big] \\ 
& = & \frac{s^*}{2} \int_{\T} \dot{\eta}^2 + 
\frac{1}{4} \int_{\T} \eta^4 
\Big[ \lambda - Q(\eta ) \Big] 
\end{eqnarray*}
where 
$$
\lambda:= \frac{\g}{2\pi |\a|} > 0 \, , \qquad 
Q(\eta ) := \frac{ \Big( \int_{\T} \eta^2 \Big)^2}{ 2 \pi \int_{\T} \eta^4} 
$$
possesses Mountain pass critical points for any $ \lambda > 0 $  
because (like in Lemma 3.14 of \cite{BB1})
$$
\inf_{\eta \in E \setminus \{0\}} Q(\eta ) = 0 \, ,  \qquad 
\sup_{\eta \in E \setminus \{0\}} Q(\eta ) = 1 
$$
(for $ \l \geq 1 $ if $ s^* = - 1 $, and 
for $ 0 < \l < 1  $ for both $ s^* = \pm 1 $). 

Such critical points satisfy 
the Euler Lagrange equation 
\begin{equation}  \label{eq di partenza}
-s^* \ddot{\eta} - \la \eta^2 \ra \eta + \lambda \eta^3 =0
\end{equation}
but their non-degeneracy is not obvious.
For this, it is convenient to express
this solutions in terms of the Jacobi elliptic sine.

\begin{proposition} \label{prop doppia}
(i) Let $s^*=-1$. Then for every $\lambda \in (0,+\infty)$ there exists an odd, analytic, $2\pi$-periodic solution $ g(t)  $ 
of (\ref{eq di partenza}) which is non-degenerate in $E $.
$ g(t) = V \sn (\Om t, m)$\, 
for   $V >0 $,\, $\Om>0$, $m \in(-\infty,-1)$ suitable constants.

\noindent
(ii) Let $s^*=1$. Then for every $\lambda \in (0,1)$ there exists 
an odd, analytic, $2\pi$-periodic solution $ g(t)  $ 
of (\ref{eq di partenza}) which is non-degenerate in $ E $. 
$ g(t) = V \sn (\Om t, m)$\, for $V >0 $,\, $\Om>0$, $m \in (0,1)$ suitable constants.
\end{proposition}

We prove Proposition \ref{prop doppia} in several steps. First we construct the solution $g$ like in Lemma \ref{jac}.

\begin{lemma} \label{jac caso u2u3}
(i) Let $s^*=-1$. Then for every $\lambda \in (0,+\infty)$ there exist
$V>0$,\, $\Om>0$, $m \in(-\infty,-1)$ such that  $g(t)=V \sn(\Omega t,m)$ is an odd, analytic, $2\pi$-periodic solution  
of (\ref{eq di partenza}) with a pole in 
$ \frac{i}{\Omega \sqrt{1-m}} \,K \big( \frac{1}{1-m} \big) $.

\noindent
(ii) Let $s^*=1$. Then for every $\lambda \in (0,1)$ there exist  
$V>0$,\, $\Om>0$, $m \in(0,1)$ such that $g(t)=V \sn(\Omega t,m)$ is an odd, analytic, $2\pi$-periodic solution  
of (\ref{eq di partenza}) with a pole in $iK(1-m)/\Omega$. 
\end{lemma}

\begin{pf}
We know that $ g_{(V,\Om,m)}(t) := V \sn (\Omega t , m )$\,
is an odd, $(4K(m)/\Om)$-periodic solution of (\ref{eq di Vsn(Omega t,m)}), see Lemma \ref{jac}. So  
it is a solution of (\ref{eq di partenza}) if $(V,\Om,m)$ verify 
\begin{eqnarray}   \label{ipotesi H caso u2u3}
\begin{cases}
  	\Omega^2 (1+m) = s^* V^2 \la \sn^2 (\cdot, m) \ra \\
   2m \Om^2 = s^* V^2 \lambda \\
   2 K(m)= \Om \pi \, .  
   \end{cases}
\end{eqnarray}
Conditions (\ref{ipotesi H caso u2u3}) give the connection between $\lambda$ and $m$:
\begin{equation}  \label{relazione lambda m}
\lambda=\frac{2m}{1+m} \, \la \sn^2 (\cdot, m) \ra  \,. 
\end{equation}
Moreover system  (\ref{ipotesi H caso u2u3}) imposes
\begin{eqnarray*} 
  \begin{cases}
	m \in (-\infty,-1) & \mathrm{if}\;\, s^*=-1 \\
	m \in (0,1) & \mathrm{if}\;\, s^*=1\,.
  \end{cases}
\end{eqnarray*}
We know that $m \mapsto  \la \sn^2 (\cdot, m) \ra$\, is continuous, strictly increasing on $(-\infty,1)$, it tends to 0 for $m \to -\infty$ and to 1 for $m \to 1$, see Lemma \ref{proprieta di sn2}. So the right-hand side of (\ref{relazione lambda m}) covers $(0,+\infty)$ for $m \in (-\infty,0)$, and it covers $(0,1)$ for $m \in (0,1)$. For this reason 
for every $\lambda >0$ there exists a unique $\bar{m} < -1$ satisfying (\ref{relazione lambda m}), and for every $\lambda \in (0,1)$ there exists a unique $
\bar{m} \in (0,1)$ satisfying (\ref{relazione lambda m}). 

The value $\bar{m}$ and system (\ref{ipotesi H caso u2u3}) determine uniquely the values $\bar{V}$, $\bar{\Omega}$.  

Analyticity and poles follow from \cite{Handbook}, 16.2, 16.10.2, pp.570,573.
\end{pf}
 
Now we have to prove the non-degeneracy of $g$. The linearized equation of (\ref{eq di partenza}) at $g$ is 
\begin{equation}  \label{linearized caso u2u3}
\ddot{h} +s^* \big( \la g^2 \ra -3\lambda g^2 \big) h = -2s^* \la gh \ra g.
\end{equation}

Let $L$ be the Green operator, i.e. for $ f \in E $, let $ H := L(f) $ be the unique solution belonging to $ E $  
of the non-homogeneous linear system
\[
{\ddot H} + s^* \big( \la g^2 \ra -3\lambda g^2 \big)\,H = f \,.
\]
We can write (\ref{linearized caso u2u3}) as
\be \label{eql caso u2u3}
 h = -2s^* \la gh \ra L(g) \, .
\ee
Multiplying by $g$ and integrating we get
\[
\la gh \ra \big[ 1+2s^* \la gL(g) \ra \big] =0\,.
\]
If \,$A_0 := 1+2s^* \la gL(g) \ra \neq 0$, then $\la gh \ra =0$, so by (\ref{eql caso u2u3}) $h=0$\, and the non-degeneracy is proved. 

\vspace{5pt}
It remains to show that $A_0 \neq 0$. As before, the key is to express  
$ L(g) $ in a suitable way.
We first look for a fundamental set of solutions of the homogeneous equation
\begin{eqnarray}  \label{HOM caso u2u3}
	\ddot{h} +s^* \big( \la g^2 \ra -3\lambda g^2 \big) h =0\,.
\end{eqnarray}

\begin{lemma}   \label{uv caso u2u3}
There exist two linearly independent solutions  of (\ref{HOM caso u2u3}), 
${\ov u} $ even, $2\pi$-periodic and ${\ov v}$ odd, not periodic, such that
$	\ub (0)=1$,\, $\dot{\ub}(0)=0$,\, $\vb (0)=0$,\, $\dot{\vb}(0)=1$,\, and 
\begin{equation} \label{v-v=rho u}
\vb (t+2\pi )-\vb (t) = \rho\, \ub (t) \quad \forall \,t
\end{equation}
for some $\rho \neq 0$.\, Moreover there hold the following expressions for $\ub$, $\vb$:
\begin{eqnarray} 
	&& \ub(t)={\dot g}(t) / {\dot g}(0)=\dot{\sn}(\Ob t, \mb) \label{eq di u}\\
	&& \vb(t)=\frac{1}{\Ob (1-\mb)}\, \sn(\Ob t)+\frac{\mb}{\mb -1}\,\dot{\sn}(\Ob t) 
	\Big[ \,t+\frac{1+\mb}{\Ob} \int_0^{\Ob t} \frac{\sn^2(\xi,\mb)}{\dn^2(\xi,\mb)}\, d\xi \Big]\,.  \label{eq di v} 
\end{eqnarray}
\end{lemma}

\begin{pf} $g$ solves (\ref{eq di partenza}) so $\dot{g}$ solves (\ref{HOM caso u2u3}); normalizing we get (\ref{eq di u}). 

By (\ref{eq di Vsn(Omega t,m)}), the function $y(t)=V \sn(\Om t,m)$ solves 
\begin{equation}   \label{aux}
\ddot{y}+s^*\la g^2 \ra y - s^* \lambda y^3=0
\end{equation}
 if $(V,\Om,m)$ satisfy 
\begin{eqnarray}  \label{per costruz di v}
	\begin{cases}
	\Om^2 (1+m) =s^* \la g^2 \ra \\
	2m \Om^2 =s^* V^2 \lambda \,.
\end{cases}
\end{eqnarray}
We solve (\ref{per costruz di v}) w.r.t. $m$ finding the one-parameter family $(y_m)$ of odd periodic solutions of (\ref{aux}), $y_m(t)=V(m) \sn (\Om(m)t,m)$. So $l(t):=(\partial_m y_m)_{|m=\mb}$ solves
(\ref{HOM caso u2u3}). We normalize $\vb(t):=l(t)/\dot{l}(0)$ and we compute the coefficients differentiating  (\ref{per costruz di v}) w.r.t. $m$. From the definitions of the Jacobi elliptic functions it holds
\[
\partial_m \sn(x,m)=-\dot{\sn}(x,m) \frac12 \int_0^x \frac{\sn^2 (\xi,m)}{\dn^2 (\xi,m)}\,d\xi\,;
\]
thanks to this formula we obtain (\ref{eq di v}).

Since $2\pi \Ob =4K(\mb)$ is the period of the Jacobi functions $\sn$ and $\dn$, by (\ref{eq di u}),(\ref{eq di v}) we obtain (\ref{v-v=rho u}) with
\[
\rho=\frac{\mb}{\mb-1}\,2\pi \Big( \,1+(1+\mb) \la \frac{\sn^2}{\dn^2} \ra \Big)\,.
\]
If $s^*=1$, then $\mb \in (0,1)$ and directly we can see that $\rho <0$. 
If $s^*=-1$, then $\mb<-1$. From the equality \,$
\la \sn^2 / \dn^2 \ra = (1-m)^{-1}\,\big(1-\la \sn^2 \ra \big) $\, 
(see \cite{Baldi}, Lemma 3, (L.2)), it results $\rho >0$. 
\end{pf}

We can note that the integral representation (\ref{Formula di L}) of the Green operator $L$ holds again in the present case. The proof is just like in Lemma  \ref{lemma della formula integrale di L}. 

\begin{lemma} \label{mostro di conti}
We can write $A_0 := 1+2s^* \la gL(g) \ra$ as function of \,$\lambda$, $\mb$,
\begin{eqnarray}  \label{eq di A0 in lambda,m}
  A_0=\frac{\lambda (1-\mb)^2 q -(1-\lambda)^2 (1+\mb)^2 +\mb q^2}{\lambda (1-\mb)^2 q}\,, 
  \qquad \quad
  q=q(\lambda,\mb):= 2-\lambda \frac{(1+\mb)^2}{2\mb} > 0\,.	
\end{eqnarray}
\end{lemma}

\begin{pf}
First, we calculate $\la gL(g) \ra$ with the integral formula (\ref{Formula di L}) of $L$. The equalities (\ref{finq}),(\ref{pe}) still hold, while similar calculations give 
\begin{eqnarray*}  
	\int_0^{2\pi} g^3 \vb =-s^* \frac{\dot{g}(0) \rho}{2 \lambda}
\end{eqnarray*}
instead of (\ref{last}). So
\begin{eqnarray}  \label{gLg intermedio}
	\la gL(g) \ra= -s^*\,\frac{\rho}{4\pi \lambda} + \frac{1}{2\pi \rho}\, \Big( \int_0^{2\pi} g \vb \Big)^2 
\end{eqnarray}
and the sign of $A_0$ is not obvious. We calculate $\int_0^{2\pi} g \vb$\, recalling that $g(t)=\Vb \sn (\Ob t,\mb)$, using formula (\ref{eq di v}) for $\vb$ and integrating by parts
\[
\int_0^{2\pi} \sn(\Ob t) \dot{\sn} (\Ob t) \mu (t) \,dt = -\frac{1}{2\Ob}\,\int_0^{2\pi} \sn^2 (\Ob t) \dot{\mu} (t) \, dt
\]
where $\mu(t):=t+(1+\mb)\Ob^{-1} \int_0^{\Ob t} \sn^2(\xi)/\dn^2(\xi)\, d\xi $.\, From \cite{Baldi}, (L.2),(L.3) in Lemma 3,  we obtain the formula
\begin{eqnarray*}  
	\la \frac{\sn^4}{\dn^2} \ra= \frac{1+(m-2) \la \sn^2 \ra}{m(1-m)}
\end{eqnarray*}
and consequently
\begin{eqnarray}
	\int_0^{2\pi} g \vb =\frac{\pi \Vb}{\Ob (1-\mb)^2}\, \big( 1+\mb -2\mb \la \sn^2 \ra \big)\,.
\end{eqnarray}
By the second equality of (\ref{ipotesi H caso u2u3}) and (\ref{gLg intermedio}) we get
\begin{eqnarray}  \label{A0 intermedio}
	A_0= 1+ \frac{2}{\l}\, \Big[ -\frac{\rho}{4\pi} + \frac{\pi \mb}{\rho (1-\mb)^4}\, \big( 1+\mb -2\mb \la \sn^2 \ra \big)^2 \Big]
\end{eqnarray}
both for $s^*=\pm 1$. From the proof of Lemma \ref{uv caso u2u3}\, we have $\rho=-2\pi \mb q\,(1-\mb)^{-2}$,\, where $q$ is defined in (\ref{eq di A0 in lambda,m}); inserting this expression of $\rho$ in (\ref{A0 intermedio}) we obtain (\ref{eq di A0 in lambda,m}).

Finally, for $\mb<-1$ we have immediately $q>0$, while for $\mb \in (0,1)$ we get $q=2-(1+\mb)\la \sn^2 \ra$  by (\ref{relazione lambda m}). Since $\la \sn^2 \ra <1$,\, it results $q>0$.  
\end{pf}
 
\begin{lemma} \label{lemma A0 neq 0}
$A_0 \neq 0 $. More precisely, sign$(A_0) = -s^*$.  
\end{lemma}

\begin{pf}
From (\ref{eq di A0 in lambda,m}), $A_0 >0$ iff $\lambda (1-\mb)^2 q -(1-\lambda)^2 (1+\mb)^2 +\mb q^2 >0$. This expression is equal to $-(1-\mb)^2 p$\,, where 
\begin{eqnarray*}  
	p=p(\lambda, \mb)=\frac{(1+\mb)^2}{4\mb}\, \lambda^2 -2 \lambda +1\,,
\end{eqnarray*}
so $A_0>0$ iff $p<0$. The polynomial $p(\l)$ has degree 2 and its determinant is 
$\Delta=-(1-\mb)^2 /\mb$. So, if $s^*=1$, then $\mb \in (0,1)$, $\Delta<0$ and $p>0$, so that $A_0<0$.

It remains the case $s^*=-1$. For $\l >0$, we have $p(\l)<0$ iff $\l > x^*$, where 
$x^*$ is the positive root of $p$,\, $x^*:=2R (1+R)^{-2}$,\, $R:= |\mb|^{1/2}$. By (\ref{relazione lambda m}), $\l>x^*$ iff 
\begin{eqnarray}  \label{6:32}
	\la \sn^2(\cdot,\mb) \ra > \frac{R-1}{(R+1)R}\,.
\end{eqnarray}
By formula (\ref{media di sn2}) and by definition of complete elliptic integrals $K$ and $E$ we can write (\ref{6:32}) as
\begin{eqnarray} \label{bomba}
	\int_0^{\pi/2}\! \Big( \frac{R-1}{(R+1)R}\,- \sin^2\! \teta \Big) \, 
	\frac{d \teta}{\sqrt{1+R^2 \sin^2 \!\teta }} \, <0 \,.
\end{eqnarray}
We put $\s:= R-1/(R+1)R$ and note that $\s<1/2$\, for every $R>0$.

$\s-\sin^2 \!\teta >0$\, iff $\teta \in (0,\teta^*)$, where $\teta^*:= \arcsin (\sqrt{\s})$, i.e. $\sin^2\! \teta^* = \s$.
Moreover $1<1+R^2 \sin^2\! \teta < 1+R^2$ for every $\teta \in (0,\pi/2)$. So
\begin{eqnarray}  \label{fine 1}
	\int_0^{\pi/2} \! \frac{\s - \sin^2\! \teta}{\sqrt{1+R^2 \sin^2\! \teta}}\,d \teta \: <
	\,\int_0^{\teta^*} \! \big(\s-\sin^2\! \teta \big)\, d \teta + \int_{\teta^*}^{\pi/2} \! \frac{\s - \sin^2\! \teta}{\sqrt{1+R^2 }}\,d \teta \,. 
\end{eqnarray}
Thanks to the formula
\[
\int_a^b \sin^2\! \teta \, d \teta \,= \frac{b-a}{2} - \frac{ \sin (2b) - \sin(2a)}{4}
\]
the right-hand side term of (\ref{fine 1}) is equal to
\[
\frac{\sin(2\teta^*)}{4}\,\Big[ (2\s-1) \Big( \frac{2\teta^*}{\sin(2\teta^*)}\,+ \frac{1}{\sqrt{1+R^2}}\,\frac{\pi-2\teta^*}{\sin(2\teta^*)} \Big) + \Big( 1-\frac{1}{\sqrt{1+R^2}} \Big) \Big]\,.
\]
Since $2\s-1 <0$\, and $\a>\sin \a$\, for every $\a >0$,\, this quantity is less than
\[
\frac{\sin(2\teta^*)}{4}\,\Big[ (2\s-1)\Big( 1+\frac{1}{\sqrt{1+R^2}} \Big) +\Big( 1-\frac{1}{\sqrt{1+R^2}} \Big)\Big]\,.
\]
By definition of $\s$, the last quantity is negative for every $R>0$, so (\ref{bomba}) is true. Consequently $\lambda >x^*$,\, $p<0$\,  and $A_0>0$. 
\end{pf}

\vspace{6pt}

As Appendix, we show the properties of the function $m \mapsto \la \sn^2(\cdot,m) \ra$ used in the proof of Lemma \ref{jac caso u2u3}.

\begin{lemma} \label{proprieta di sn2}
The function $\ph:(-\infty,1) \rightarrow \R$\,,\, $m \mapsto \la \sn^2(\cdot,m) \ra$\, is continuous, differentiable, strictly increasing, and \,$\lim_{m \to -\infty} \ph(m)=0$\,,\, $\lim_{m \to 1} \ph(m)=1$\,. 
\end{lemma}

\begin{pf}
By (\ref{media di sn2}) and by definition of complete elliptic integrals $K$ and $E$,
\begin{eqnarray*}  
	\ph(m)=\frac{K(m)-E(m)}{mK(m)}= \,\int_0^{\pi/2} \!\!\!\! \frac{\sin^2 \! \teta \, d\teta}{\sqrt{1-m\sin^2 \!\teta}}\, \: \Big( \int_0^{\pi/2} \!\!\!\! \frac{d\teta}{\sqrt{1-m\sin^2 \!\teta}}\, \Big)^{-1}\,, 
\end{eqnarray*}
so the continuity of $\ph$ is evident. 

Using the equality $\sin^2\!+\cos^2\!=1$\, and the change of variable $\teta \rightarrow \pi/2 \, - \teta$\, in the integrals which define $K$ and $E$, we obtain the formulae
\begin{eqnarray}  \label{formula magica foglio verde}
	K(m)=\frac{1}{\sqrt{1-m}}\,K\Big( \frac{m}{m-1} \Big)\,, \qquad \quad 
	E(m)= \sqrt{1-m}\, E \Big( \frac{m}{m-1} \Big)\, \qquad \forall \, m <1\,.
\end{eqnarray}
We put $\mu:=m/(m-1)$, so it results
\begin{eqnarray}  \label{con mu}
	\ph(m)=1-\frac{1}{\mu}+\frac{E(\mu)}{\mu K(\mu)}\,.
\end{eqnarray}
Since $\mu$ tends to 1 as $m \to -\infty$,\, $E(1)=1$\, and \,$\lim_{\mu \to 1} K(\mu)= +\infty$,\, (\ref{formula magica foglio verde}),(\ref{con mu}) give \,$\lim_{m \to -\infty} \ph(m)=0$. 
Since $E(m)/K(m)$\, tends to 0 as $m \rightarrow 1$, (\ref{media di sn2}) gives \,$\lim_{m \to 1} \ph(m)=1$.

Differentiating the integrals which define $K$ and $E$ w.r.t. $m$ we obtain the formulae
\begin{eqnarray*}  
	E'(m)=\frac{E(m)-K(m)}{2m}\,, \qquad \quad 
	K'(m)=\frac{1}{2m}\, \Big( \int_0^{\pi/2} \!\!\!\! \frac{d \teta}{(1-m\sin^2\!\teta )^{3/2}} \,-K(m) \Big)\,,
\end{eqnarray*}
so the derivative is
\[
\ph'(m)=\frac{1}{2m^2 K^2(m)}\, \Big[ \,E(m) \int_0^{\pi/2} \!\!\!\! \frac{d \teta}{(1-m\sin^2\!\teta )^{3/2}}\,- K^2(m)\, \Big]\,.
\]
The term in the square brackets is positive by strict H$\ddot{\mathrm{o}}$lder  inequality for $(1-m\sin^2\!\teta)^{-3/4}$ and $(1-m\sin^2\!\teta)^{1/4}$.
\end{pf}

\vspace{6pt}

\noindent
{\bf Acknowledgements:}
The authors thank Philippe Bolle for useful comments.

\end{document}